\documentclass[10pt]{article}

\usepackage{amsmath,amssymb,amsfonts, slashed, shuffle}
\usepackage{dsfont}

\newtheorem{theorem}{Theorem}
\newtheorem{proposition}{Proposition}
\newtheorem{definition}{Definition}
\newtheorem{corollary}{Corollary}

\title{Leading log expansion of system of combinatorial Dyson Schwinger equations}
\author{Lucas  Delage}
\date{\today}

\begin{document}

\maketitle

\begin{abstract}
We study combinatorial Dyson Schwinger equations, expressed in the Hopf algebra of words with a quasi shuffle product. We map them into an algebra of polynomials in one indeterminate $L$ and show that the leading log expansion one obtains with such a mapping are simple power law like expression.
\end{abstract}

\textbf{Acknowledgment}

\textit{I want to thanks Dirk Kreimer and his groups for their help and their friendship which made my work with them incredibly nice}

\section{Introduction}
Dyson Schwinger equations are functional equations in Quantum Field Theory fulfilled by the Green's functions, which are the expectation values of fields monomials and lead to some scattering amplitudes in particle physics. These equations are fix-point and self coherent equations. They also have analogies in statistical field theories. Dirk Kreimer, in his works about Hopf algebraic renormalization, extended these Dyson Schwinger equations into a combinatorial form, using a Hochschild-1-cohomology, often called $B_+$ operator, which plays the role of the linear forms in the functional expression. The equations become equations for series $X^r$'s expressed in a specified Hopf algebra. It leads to rich purely mathematical problems, which were for example studied by Loic Foissy in [6] or more recently by Joachim Kock in [5]. 

Here, we study a specific system of combinatorial Dyson Schwinger equations in the Hopf algebra of words endowed with a quasi shuffle product, in order to find at the end the leading series of the corresponding Green's functions, called leading log expansion. The Green's functions are supposed to depend of only one kinematical variable $L$, and one coupling constant $\alpha$. This simplifies the correspondance, given by the so-called Feynman rules, between combinatorics and physics. Given these settings, we write our CDSE's as:
\begin{eqnarray}
\nonumber X^r = \mathbf{1}_W + sign(\eta_r)\alpha B_+^{a_r}(X^r Q), \\
\nonumber Q = \prod_{r'=1}^R (X^{r'}) ^{\eta_{r'}}.
\end{eqnarray}
The work done here is largely inspired by [7]. The reason of why we use this Hopf algebra of words is because it allows a nice factorization of words in components only linear in L once mapped by the Feynman rules.

We first start by giving the definitions and the elementary properties we will need. They are reminders on the theory of Hopf algebras and Hopf algebraic renormalization, developed in [3], [4], [8] and [9]. We then define what is the Hopf algebra of words that we will use to define our CDSE's, and finally pick out their leading log expansion $G_{LL}^r$ thanks to a theorem proved in the Appendix. We finally  show that: 
\begin{eqnarray}
\nonumber G_{LL}^r &=& (1 + A \alpha L)^{-\frac{c_r}{A}} \text{  if } \eta_r < 0, \\
\nonumber G_{LL}^r &=& (1 - A \alpha L)^{\frac{c_r}{A}} \text{  if } \eta_r > 0 ,\\
\nonumber A &=& \sum_{r'=1}^R \eta_{r'} c_{r'}.
\end{eqnarray}

\section{Preliminaries}

We recall there some useful definitions and properties about Hopf algebras. If the reader is not familiar with the concept of Hopf algebra, [8] may help him a lot. It is also in this paper that one can find the proofs of the next properties.

\begin{definition}
A bialgebra $(H,m,\mathbf{1},\Delta,\hat{\mathbf{1}})$ is a vector space $H$ over a field $\mathbb{K}$ together with an associative product $m : H \otimes H \rightarrow H$, a unit $\mathbf{1}$ such that
$\forall x \in H, m(x \otimes \mathbf{1}) = m(\mathbf{1} \otimes x) = x$, a coassociative coproduct $\Delta : H \rightarrow H \otimes H$ and a counit $\hat{\mathbf{1}}$ satisfaying $(Id \otimes \hat{\mathbf{1}}) \circ \Delta (x) = (\hat{\mathbf{1}} \otimes Id) \circ \Delta (x) = x$. Furthermore , $\Delta$ and $\hat{\mathbf{1}}$ have to be algebra morphisms with respect to the product $m$ or equivalently, $m$ and $\mathbf{1}$ have to be coalgebra morphisms with respect to the coproduct $\Delta$
\end{definition}

\begin{definition}
If $(H,m,\mathbf{1},\Delta,\hat{\mathbf{1}},S)$ is a bialgebra together with a linear map $S : H \rightarrow H$ which fulfills $m(S \otimes Id)\circ \Delta = m(Id \otimes S) \circ \Delta = \mathbf{1} \hat{\mathbf{1}}$, it is called a Hopf algebra. S is called the antipode.
\end{definition}

\begin{definition}
A bialgebra (respectively, a Hopf algebra) is called filtered iff there exist subspaces $H^0 \subset H^1 \subset ... \subset H^n \subset...$ such that $\bigcup_n H^n = H$, $m(H^p \otimes H^q) \subset H^{p+q}$, $\Delta(H^n) \subset \sum_{p+q=n} H^p \otimes H^q$ (respectively, if furthermore $S(H^n) \subset H^n)$. 
\end{definition}

\begin{definition}
A bialgebra (respectively, a Hopf algebra) is called graded iff there exist subspaces $H_0, H_1, ..., H_n, ...$ such that $\oplus_n H_n = H$, $m(H_p \otimes H_q) \subset H_{p+q}$, $\Delta(H_n) \subset \oplus_{p+q=n} H_p \otimes H_q$ (respectively, if furthermore $S(H_n) \subset H_n)$. 
\end{definition}

\begin{definition}
A filtered (respectively graded) bialgebra or Hopf algebra is called connected iff $H^0$ (respectively $H_0$) is one-dimensional.
\end{definition}

\textbf{Remark:}
A graded bialgebra (respectively Hopf algebra) is in particular a filtered bialgebra (respectively Hopf algebra). The canonical filtration associated with a grading is given by, keeping the same notation as above:
\begin{equation}
H^n = \oplus_{j=0}^n H_j.
\end{equation}

\begin{definition}
A bialgebra or a Hopf algebra is called pointed iff all its simple left (or right) comodules are one-dimensional.
\end{definition}

\begin{definition}
Let $(A, m_A, \mathbf{1}_A)$ be a unital algebra and $H$ be a bialgebra or a Hopf algebra as above. One defines the convolution product $\star$ on $Hom(H,A)$ by:
\begin{equation}
\forall f_1, f_2 \in Hom(H,A), f_1 \star f_2 = m_A(f_1 \otimes f_2) \circ \Delta .
\end{equation}
\end{definition}

\textbf{Remark:}
The convolution product admit a neutral element: $\mathbf{1}_A \hat{\mathbf{1}}$.

\vspace{10pt}
\textbf{Remark:}
The antipode of a Hopf algebra is the inverse of the identity for the convolution product in $Hom(H,H)$.

\begin{definition}
An element $x \in H$ such that $\Delta(x) = \mathbf{1} \otimes x + x \otimes \mathbf{1}$ is called a primitive element. The subspace of $H$ of primitive elements is denoted by $Prim(H)$.
\end{definition}

\begin{definition}
The reduced coproduct $\widetilde{\Delta}$ is defined by $\widetilde{\Delta}(\mathbf{1}) = 0$ and by:
$$ \forall x \neq \mathbf{1}, \hspace{2pt} \widetilde{\Delta}(x) = \Delta(x) - \mathbf{1} \otimes x - x \otimes \mathbf{1}, $$
and the $k^{th}$ power of the reduced coproduct by:
$$\widetilde{\Delta}^k = (\underbrace{Id \otimes Id \otimes ... \otimes Id }_{k-1 \hspace{2pt} times}\otimes \widetilde{\Delta})...(Id \otimes \widetilde{\Delta})\widetilde{\Delta}. $$
\end{definition}

One can easily check that the reduced coproduct is coassociative.

\begin{proposition}
If H is a pointed bialgebra or a pointed Hopf algebra, one can endowed it with the coradical filtration:
$$x \in H^k \text{ iff } \widetilde{\Delta}^k(x) = 0.$$
\end{proposition}

An element $x$ such that $\widetilde{\Delta}^k(x) = 0$ and $\widetilde{\Delta}^{k-1} \neq 0$ is said to be of coradical degree $k$.

\begin{proposition}
If H is a connected filtered bialgebra, then it extends canonically to a Hopf algebra. The antipode is defined by $S= \sum_k (\mathbf{1} \hat{\mathbf{1}} - Id)^{\star k}$. It is given by $S(\mathbf{1}) =\mathbf{1}$ and for $x \neq \mathbf{1}$ recursively by one of the two following formulas:
\begin{eqnarray}
S(x) &=& -x -(S \otimes Id)\tilde{\Delta}(x), \\
S(x) &=& -x -(Id \otimes S)\tilde{\Delta}(x).
\end{eqnarray} 
\end{proposition}

\subsection{Hochschild Cohomology}

The work we present here can also be found in [9] or in [4].

One defines a Hochschild cochain complex:
\begin{equation}
\lbrace Hom(H, H^{\otimes n}) \rbrace _{n \in \mathbb{N}}
\end{equation}
with coboundary maps $b_n$ such that if we let $L_n \in Hom(H,H^{\otimes n})$ and $h\in H$ then:
\begin{eqnarray}
b_n L_n (h) = (id \otimes L_n) \Delta (h) + \sum_{i=1}^n (-1)^i \Delta_i L_n (h) +(-1)^{n+1} L(h) \otimes \mathbf{1}, \\
\Delta_i = \underbrace{(id \otimes .... \otimes id \otimes \Delta \otimes id \otimes ... \otimes id)}_{n \text{  times}} \text{with $\Delta$ at the $i^{th}$-place.} 
\end{eqnarray}
The elements of $\ker(b_{n+1})$ are called Hochschild-n-cocycles and the set of all such Hochschild-n-cocycles is denoted by $ZH^{n}$. The elements of $Im(b_n)$ are called Hochschild-n-coboundaries and the set of all such Hochschild-n-coboundaries is denoted by $BH^n$. Finally, the $n^{th}$ cohomology space is $HH^n = \frac{ZH^n}{BH^n}$.

\vspace{10pt}
In this paper, we are interested by Hochschild-1-cocycles, which we will denote by $B_+^x$ such that $B_+^x (\mathbf{1}) = x \in Prim(H) $. As we have for any $L_0 \in Hom(H,\mathbb{K})$
\begin{equation}
b_0 \circ L_0 (\mathbf{1}) = L_0(\mathbf{1})\mathbf{1} - (id \otimes L_0)\circ \Delta(\mathbf{1}) = 0
\end{equation}
our $B_+^x$ is in fact a Hochschild-1-cocycle. We also have \begin{equation}
\Delta \circ B_+(\mathbf{1}) = B_+(\mathbf{1}) \otimes \mathbf{1} + \mathbf{1} \otimes B_+(\mathbf{1}),
\end{equation}
i.e. $B_+(\mathbf{1}) \in Prim(H)$ for all Hochschild-1-cocycles $B_+$.

\subsection{Group of characters}

Let $(A, m_A, \mathbf{1}_A)$ be a unital commutative algebra, $(H, m, \mathbf{1}, \Delta, \hat{\mathbf{1}}, S)$ be a Hopf algebra, and $\star$ be the convolution product as before. Consider the set $$G(H,A) = \lbrace \phi \in Hom(H,A) / \phi(\mathbf{1}) = \mathbf{1}_A \rbrace .$$ Elements of $G(H,A)$ are called characters.

\begin{proposition}
The set $(G(H,A), \star, \mathbf{1}_A \hat{\mathbf{1}})$ is a group 
\end{proposition}

Furthermore, consider the set $$g(H,A) = \lbrace \sigma \in Hom(H,A) / \sigma(xy) = \sigma(x) \hat{\mathbf{1}}(y) + \hat{\mathbf{1}}(x) \sigma(y) \rbrace$$ and the bracket $[,]$ with $$ \forall \sigma, \rho \in Hom(H,A), [\sigma,\rho] = \sigma \star \rho - \rho \star \sigma.$$ 

\begin{proposition}
The set $(g(H,A, [,])$ is a Lie Algebra.
\end{proposition}
Elements of $g(H,A)$ are called infinitesimal characters.

\vspace{10pt}
\textbf{Remark: } $\forall \sigma \in g(H,A), \sigma(\mathbf{1}) = 0$

\begin{proposition}
\begin{eqnarray}
\forall \sigma \in g(H,A), \exp_{\star}(\sigma) = \sum_n \frac{\sigma^{\star n}}{n!} \in G(H,A). \\
\forall \phi \in G(H,A), \exists \sigma \in g(H,A) / \phi = \exp_{\star}(\sigma).
\end{eqnarray}
\end{proposition}
The sum in (10) is bounded by the coradical degree of $x \in H$ as $\sigma(\mathbf{1}) = 0$. The proofs of all these propositions can be found in [4] and in [8].

\section{Hopf algebra of words}

Let $\Omega$ be a countable space. We define a symmetric function: $\Omega \otimes \Omega \rightarrow \Omega$. We call inherited elements those which are in $Im(\Theta) \subset \Omega$. We define an equivalence relation: $$E: \Theta(a,\Theta(b,c)) \sim \Theta(\Theta(a,b),c)) .$$ We let $H_{\Theta}$ be the quotient space $\Omega / E$ in which $\Theta$ is fully symmetric and in this space we write $\Theta(\Theta(a,b),c) \equiv \Theta(a,b,c)$. Finally, we suppose that there are $R \in \mathbb{N}^*$ fixed non inherited elements in $H_{\Theta}$. The set of these non inherited elements is denoted by $H_L$.

We call the set $H_{\Theta}$ an alphabet and any of its element a letter. We call a word a concatenation of several letters and we denote by $H_W$ the vector space generated by all words; each word defining a basis element.
We say that a word is of length $l(w) \in \mathbb{N}$ if it's written as the concatenation of $l(w)$ letters. The length of a sum of words is defined to be the length of the word with the biggest length in the sum. We let $\varnothing$, the empty word, beeing the unique word of length 0. One define a grading on $H_W$, denoted by $\mid w\mid$, for any word $w$ by:
\begin{eqnarray}
&& \mid \varnothing \mid = 0 ,\\
\forall a \in H_L ,&& \mid a \mid = 1, \\
\forall a,b \in H_{\Theta} ,&& \mid \Theta (a,b) \mid =  \mid a \mid + \mid b \mid, \\
\forall u,v \in H_W ,&& \mid uv \mid = \mid u \mid + \mid v \mid.
\end{eqnarray}

We define recursively a map $\shuffle$, called the shuffle product, by: 
\begin{eqnarray}
\shuffle &:&  H_W \otimes H_W \rightarrow H_W  \\
\forall w \in H_W, w \shuffle \varnothing &=& w, \\ 
\forall a_i u, a_j v \in H_W, a_iu \shuffle a_jv &=& a_i(u \shuffle a_jv) + a_j(a_iu \shuffle v).
\end{eqnarray}

Similarly, we define the quasi shuffle product, denoted by $\shuffle_{\Theta}$,
\begin{eqnarray}
\shuffle_{\Theta} &:&  H_W \otimes H_W \rightarrow H_W  \\
\forall w \in H_W, w \shuffle_{\Theta} \varnothing &=& w, \\ 
\nonumber \forall a_i u, a_j v \in H_W, a_iu \shuffle_{\Theta} a_jv &=& a_i(u \shuffle_{\Theta} a_jv) + a_j(a_iu \shuffle_{\Theta} v) \\
&& + \Theta(a_i,a_j)u\shuffle_{\Theta} v.
\end{eqnarray}

We let 
\begin{eqnarray}
\nonumber \Delta &:& H_W \rightarrow H_W \otimes H_W \\
& & w \mapsto \sum_{uv = w} v \otimes u
\end{eqnarray}
be the deconcatenation. Finally, we denote by $\delta_{\varnothing}$ the indicatrix of $\varnothing$.

We give an order on $H_L$ by $a_i < a_j \Leftrightarrow i < j $ for $i,j \in \lbrace 1,...,R \rbrace$.  The order is extended on $H_{\Theta}$ by:
\begin{eqnarray}
\forall a,b,c \in H_R, && a < \Theta (b,c), \\
\forall a,b,c \in H_{\Theta}, && \Theta(a,b) < \Theta(a,c) \Leftrightarrow b<c, \\
\forall a,b,c,d \in H_{\Theta}, && \Theta(a,b) < \Theta(c,d) \Leftrightarrow a<c. 
\end{eqnarray}
We then extend this order to the words by saying $u <uv$ and $uav < ubv'$ with $u,v,v'$ being words and $a,b$ letters such that $a<b$. This order is often call the lexicographic order.

\begin{definition}
A word $w \in H_W$ is called Lyndon iff $\forall u,v \in H_W, w=uv \Rightarrow w<v$.
\end{definition}

\begin{proposition}
$(H_W, \shuffle, \varnothing, \Delta, \delta_{\varnothing})$ and $(H_W, \shuffle_{\Theta}, \varnothing, \Delta, \delta_{\varnothing})$ are bialgebras, graded connected by $\mid . \mid$. As algebras, They are freely generated by the Lyndon words.
\end{proposition}

This is proved in [1], [2].
  
\begin{corollary}
$(H_W, \shuffle, \varnothing, \Delta, \delta_{\varnothing}, S)$ and $(H_W, \shuffle, \varnothing, \Delta, \delta_{\varnothing}, S_{\Theta})$ are Hopf algebras, with S and $S_{\Theta}$ defined as in Proposition (2).
\end{corollary}

\textbf{Remark:} By definition, the shuffle product and the quasi shuffle product are commutative. They conserve the graduation $\mid . \mid$. Furthermore, the length of words defines another graduation for $(H_W, \shuffle, \varnothing, \Delta, \delta_{\varnothing})$, and only a filtration for $(H_W, \shuffle, \varnothing, \Delta, \delta_{\varnothing}, S_{\Theta})$.

\vspace{10pt}

We denote by $C(n)$ a partition of $\lbrace 1,...,n \rbrace, n \in \mathbb{N}$. We define an action of $C(n)$ over a word $w = a_{i_1}...a_{i_n}$ of length $n$ by:

\begin{equation}
<(j_1,...,j_l) \mid w> = \Theta(a_{i_1},...,a_{j_1})...\Theta(a_{j_{l+1}},...,a_{i_n}), \hspace{10pt} \forall j_1 + ... + j_l = n
\end{equation}
with the convention that $\Theta(a) = a, \hspace{5pt} \forall a \in H_L$. Here the $...$ denotes the conatenation of the several $\Theta(.)$ so obtained.

Now we let $\hat{\exp} : H_W \rightarrow H_W$ be a linear map such that $\hat{\exp}(1) = 1$ and for any non empty word $w$,

\begin{equation}
\hat{\exp}(w) \sum_{(j_1,...j_l)\in C(l(w))} \frac{1}{j_1 ! ... j_l !} <(j_1,...,j_l) \mid w> .
\end{equation}

\begin{proposition}
The map $\hat{\exp}$ is an isomorphism between $(H_W, \shuffle, \varnothing, \Delta, \delta_{\varnothing}, S)$ and $(H_W, \shuffle_{\Theta}, \varnothing, \Delta, \delta_{\varnothing}, S_{\Theta})$.
\end{proposition}

This is shown in [2].

\vspace{10pt}

The primitives elements for both Hopf algebras (as they share the same coproduct) are the letters ([1]). We will spend the next lines to caracterize some indecomposable elements. 

\begin{definition}
$\forall$ u,v $\in H_W$, we call Lie bracket the bracket $[u,v] \equiv uv - vu$.
\end{definition}

\begin{definition}
The elements generated by Lie brackets of letters are called Lie polynomials. A Lie polynomial is said to be of degree n if it is an iterated bracket of n different letters. 
\end{definition}

In particular, the Lie polynomial of degree $1$ are the letters, and if $P$ is a Lie polynomial of degree $n$, $l(P) = n$. In the following $P_n$ denotes a Lie polynomial of degree $n$.

\begin{proposition}
The Lie polynomials can be seen as the indecomposable elements of $(H_W, \shuffle, \varnothing, \Delta, \delta_{\varnothing}, S)$.
\end{proposition}

And this, instead of the Lyndon words. It will be the case in the rest of the article. We have the following corollary:

\begin{corollary}
Let P be a Lie polynomial. Then $\hat{\exp}(P)$ can be seen as an indecomposable element of $(H_W, \shuffle_{\Theta}, \varnothing, \Delta, \delta_{\varnothing}, S_{\Theta})$.
\end{corollary}

This means exactly that any word in $(H_W, \shuffle_{\Theta}, \varnothing, \Delta, \delta_{\varnothing}, S_{\Theta})$ can be written  as a sum of products of the $\hat{\exp}$ of Lie polynomials. Again, this is shown in [1] for the proposition and [2] for the corollary.

\vspace{10pt}

We ask the reader to focus his attention on a few things coming from the preceding propositions. First, the map $\hat{\exp}$ respects the graduation defining by $\mid . \mid$ as well as the coradical filtration and the filtration induced by the length of words. Secondly, if $P$ is a Lie polynomial of degree $n$ then $\mid P \mid$ is equal to the sum of the degrees (with respect to $\mid . \mid$) of the letters composing it (e.g $P = [a_{i_1},[...,[a_{i_{n-1}}, a_{i_n}]...]] \Rightarrow \mid P \mid = \mid a_{i_1} \mid + ... + \mid a_{i_n} \mid$). Furthermore, we have that $\hat{\exp}(P) = P +$ some words of length strictly less than $l(P)$. Finally, for any word $w$, $cor(w) \leq l(w) \leq \mid w \mid$ with $cor(w)$ denoting the coradical degree of $w$.

These remarks allow us to decompose any word $w$ in $(H_W, \shuffle_{\Theta}, \varnothing, \Delta, \delta_{\varnothing}, S_{\Theta})$ as follows:
\begin{equation}
w = \sum_{(i_1,..,i_l)\in C(l(w))} \sum_{j_{i_1},...,j_{i_l}} P_{j_{i_1}} \shuffle ... \shuffle P_{j_{i_l}} + \text{ words of length } < l(w).
\end{equation}
Notice carefully that we use the shuffle product in the preceding equation and not the quasi shuffle product, as the $\Theta$-part of the quasi shuffle product decreases the length of words at least by one.

Last but not least, we want to point out that for any letters $a_{i_1},...,a_{i_n}$, $a_{i_1} \shuffle... \shuffle a_{i_n}$ gives rise to a sum of $n!$ words, each of these corresponding to a permutation of $n$ elements. In other words, if we let $S_n$ be the group of permutations of $n$ elements and $s$ a permutation of these $n$ elements, one has:
\begin{equation}
a_{i_1} \shuffle... \shuffle a_{i_n} = \sum_{s \in S_n} a_{s(i_1)}...a_{s(i_n)}.
\end{equation}

\subsection{Feynman rules in the Hopf algebra of words}

In this section (and in general in all the next sections) we denote by $H_W$ the Hopf algebra  $(H_W, \shuffle_{\Theta}, \varnothing, \Delta, \delta_{\varnothing}, S_{\Theta})$; and we keep the notations of setion 1. In particular, we will denote $\varnothing$ by $\mathbf{1}$, $\delta_{\varnothing}$ by $\hat{\mathbf{1}}$ and the reduced coproduct by $\widetilde{\Delta}$.

Consider the set $\mathbb{C}[L]$ of polynomial in indeterminate L. An element $\phi \in G(H_W,\mathbb{C}[L])$ is called Feynman rules. It is generated by an element $\sigma \in g(H_W,\mathbb{C})$ by the formula 
\begin{equation}
\phi = \exp_{\star}(L \sigma).
\end{equation}
Notice that $\phi$ is entirely determined by the value of $\sigma$ on the indecomposable elements of $H_W$, namely the images by $\hat{\exp}$ of Lie polynomials. In the following, we suppose $\phi$ and $\sigma$ fixed, and suppose that $\sigma$ does not vanish on any element of $H_L$, so let say that $\forall a_i \in H_L, \sigma(a_i) = c_i \in \mathbb{C}$.

According to the Feynman rules, one can define a graduation on the Hopf algebra $H_W$, called kinematical graduation. An element $x \in H$ is called of kinematical degree $n$ if $\deg(\phi(x)) = n$, where $\deg$ denotes the usual degree of a polynomial in its indeterminate.
\begin{proposition}
An element of coradical degree $n$ is at most of kinematical degree $n$.
\end{proposition}

\textbf{Proof} :

 We have $\sigma(\mathbf{1}) = 0$, it follows:
 $$\forall x \in H, \sigma^{\star n}(x) = m(\sigma \otimes ... \otimes \sigma) \Delta^n (x) = m(\sigma \otimes ... \otimes \sigma) \widetilde{\Delta}^n (x).$$
$\square$

\begin{theorem}
Let $w = a_{i_1}...a_{i_n}$ be a connected word of length n. Then its kinematical degree is at most n and the coefficient of $L^n$ in $\phi(w)$ is proportional to $\phi(a_{i_1} \shuffle_{\Theta}... \shuffle_{\Theta} a_{i_n}) = c_{i_1}...c_{i_n}$.
\end{theorem}

\textbf{Proof} :

The first part of the theorem is trivial, as $cor(w) \leq l(w)$ and by applying proposition 9. For the second part, we decompose $w$ as in equation (29). In this decomposition, the words of length strictly less than $l(w) = n$, as they also have a coradical degree strictly less than $n$, will not contribute to the coefficient of $L^n$. Furthermore, as the target algebra $\mathbb{C}[L]$ is commutative, each Lie polynomial of degree $\geq 2$ will be mapped to 0 by $\sigma^{\star n}$. The result follows. $\square$

\subsection{Combinatorial Dyson Schwinger equations}

Before defining what are the Combinatorial Dyson Schwinger equations in the Hopf algebra of words, we have to look at what are exactly the Hochschild-1-cocycles in the Hopf algebra $(H_W, \shuffle_{\Theta}, \mathbf{1}, \Delta, \hat{\mathbf{1}}, S_{\Theta})$. Remember from 1.1 that we focus on operator $B_+^{a_r}$ such that $B_+^{a_r}(\mathbf{1}) = a_r$. Remember also that by definition, $\Delta(w) = \sum_{uv = w} v \otimes u$. So we have, for any non empty word $w$, 
\begin{eqnarray}
 \Delta B_+^{a_r}(w) &=& (id \otimes B_+^{a_r}) \Delta (w) + B_+^{a_r}(w) \otimes \mathbf{1} \\
 &=& \sum_{uv = w} v \otimes B_+^{a_r}(u)  + B_+^{a_r}(w) \otimes \mathbf{1}
\end{eqnarray}
and, in the second hand,
\begin{equation}
 \Delta B_+^{a_r}(w) = \sum_{u'v' = B_+^{a_r}(w)} v' \otimes u'
\end{equation}
In particular, one gets that $w \otimes a_r$ must appear in $ \Delta B_+^{a_r}(w)$, thus that $a_r w = B_+^{a_r}(w)$ and thus that $B_+^{a_r}$ is the operator which adds $a_r$ at the beginning of each word.

\vspace{10pt}

Let $H_W[[\alpha]]$ be the ring of formal series with parameter $\alpha$ and coefficients in $H_W$. We define Combinatorial Dyson Schwinger equations as equations for these formal series living in $H_W[[\alpha]]$. As we consider a system of $R \in \mathbb{N}^*$ such equations, we define several series denoted by  $X^r$, $r \in \lbrace 1,...,R \rbrace$. The Combinatorial Dyson Schwinger equations of our interest are:
\begin{eqnarray}
X^r = \mathbf{1}_W + sign(\eta_r)\alpha B_+^{a_r}(X^r Q), \\
Q = \prod_{r'=1}^R (X^{r'}) ^{\eta_{r'}}.
\end{eqnarray}

Thus, we extend the Feynman rules define in 2.1 to the $X \in H_W[[\alpha]]$:
\begin{eqnarray}
\nonumber \phi &:& H_W[[\alpha]] \rightarrow \mathbb{C}[[L]][[\alpha]] \\
&& X = \sum_n w_n \alpha^n \mapsto \sum_n  \phi(w_n) \alpha^n.
\end{eqnarray}
We call $G^r(\alpha,L) = \phi(X^r) = \sum_{i,j=0}^{\infty} b^r_{i,j} \alpha^i L^j$ with coefficients $b_{i,j} \in \mathbb{C}$. $G^r(\alpha,L)$ is called the log expansion of the series $X^r$.

\begin{proposition}
Consider $\lbrace X^r = \sum_n w_n^r \alpha^n \rbrace_{r \in R}$ a solution of the system of CDSE's (34). Then any of the $w_n^r$ is of length $n$.
\end{proposition}

\textbf{Proof:} 

We proceed by induction on the order in $n$ of $X^r$. For $n=0$, $w_0^r = \mathbf{1}_W$ for any $r$. Now consider that for a fixed $n$ any of the $w_m^r, m \leq n$ is of length $m$. Now $w_{n+1}^r$ is the coefficient of the term $\alpha^{n+1}$ in $X^r$. Thus it is equal to $\alpha B_+^{a_r}(W_n^r)$ where $W_n^r$ is the term of order $\alpha^n$ in the product $X^{r}Q$. $X^{r}Q$ can be expanded in a formal serie in each of its variable, so we can write $W_n^r = \sum_i \prod_{j,r'} w_{i,j}^{r'}$ with $\sum j = n$, respecting the order in $\alpha$. Hence $W_n^r$ is of length $n$ for each $r$. As $B_+^{a_r}$ increases the length by 1, $w_{n+1}^r$ is of length $n+1$ and the proposition is proved. $\square$

\vspace{10pt}
\textbf{Remark:} By the same induction procedure, one can easily show that any of the $w_n^r, n \geq 1$ can be written as a word beginning by $a_r$.

\vspace{10pt}
The following is a direct consequence of the preceding proposition and theorem 1.

\begin{corollary}
Let  $\lbrace X^r = \sum_n w_n^r \alpha^n \rbrace_{r \in R}$ a solution of the system of CDSE's (34). Then any $G^r(\alpha,L)$ can be written as $\sum_{i=0}^{\infty} \sum_{j=0}^i b^r_{i,j} \alpha^i L^j$.
\end{corollary}
We call $G^r_{LL}(\alpha,L)$  the sum $\sum_{i=0}^{\infty} b^r_{i,i} \alpha^i L^i$. $G^r_{LL}(\alpha,L)$ is called the leading log expansion of the series $X^r$.

\begin{proposition}
Let w = $a_{i_1}\shuffle_{\Theta}...\shuffle_{\Theta} a_{i_n}$, $a_{i_1},..,a_{i_n} \in H_L$. Then the coefficient of the term in $L^{n+1}$ in $\phi(B_+^{a_r}(w))$ is equal to $\frac{1}{n+1}\phi(a_{i_1} \shuffle_{\Theta}... \shuffle_{\Theta} a_{i_n})$.
\end{proposition}

\textbf{Proof:}

Denoting as above by $S_n$ the group of permutations of $n$ elements, we know that $a_{i_1}\shuffle_{\Theta}...\shuffle_{\Theta} a_{i_n} = \sum_{s \in S_n} a_{s(i_1)}...a_{s(i_n)}$ + words of length at most $n-1$. By applying $B_+^{a_r}$ on these words, we obtain $n!$ words beginning with $a_r$ of length $n+1$. 

Now consider the product $a_r \shuffle_{\Theta} a_{i_1}\shuffle_{\Theta}...\shuffle_{\Theta} a_{i_n}$. You obtain $n!$ words beginning with $a_r$, $n*n!$ other words of length $n+1$ and some other words of length at most $n$, which will be mapped to polynomial of degree at most $n$ by $\phi$. For the other words of length $n+1$, bring $a_r$ to the first place by applying the transformation $a_ia_r = a_ra_i + [a_i,a_r]$. You then obtain $(n+1)!$ words beginning with $a_r$, and $n*n!$ words containing Lie polynomials of degree $n+1$ which will be mapped again to a polynomial of degree $n$ in $L$ by $\phi$ thanks to the same argument as in the proof of Theorem 1. In fact, the $(n+1)!$ words beginning with $a_r$ are exactly $n+1$ times the words we obtained by computing $B_+^{a_r}(w)$. Applying Theorem 1 on both sides ends to prove the result. $\square$

\vspace{10pt}

\begin{theorem}
Consider $\lbrace X^r = \sum_n w_n^r \alpha^n \rbrace_{r \in R}$ a solution of the system of CDSE's (34). Then the leading log expansion of each $X^r$ obeys
\begin{eqnarray}
\nonumber \sum_{q_1,...,q_R ; q_r \geq 1} (q_1+...+q_R) C^r_{q_1,...,q_R} \phi(a_1)^{q_1}...\phi(a_r)^{q_r -1}...\phi(a_R)^{q_R} \alpha^{q_1+...+q_R} \\
 = G_{LL}^r \prod_{r'=1}^R (G_{LL}^{r'})^{\eta_{r'}}
\end{eqnarray}
with some coefficients $C^r_{q_1,...,q_R} \in \mathbb{C}$.
\end{theorem}

This is proved in a more general case in the appendix.

\vspace{10pt}

Going on with this relation, one gets :
\begin{eqnarray}
\nonumber \sum_{q_1,...,q_R ; q_r \geq 1} (q_1+...+q_R) C^r_{q_1,...,q_R} (c_1L\alpha)^{q_1}...(c_rL\alpha)^{q_r -1}...(c_R L\alpha)^{q_R} \\
= G_{LL}^r \prod_{r'=1}^R (G_{LL}^{r'})^{\eta_{r'}}. 
\end{eqnarray}
The left-hand-side can be written as
\begin{eqnarray}
&& \sum_{q_1,...,q_R ; q_r \geq 1} (q_1+...+q_R) C^r_{q_1,...,q_R} (c_1L\alpha)^{q_1}...(c_rL\alpha)^{q_r -1}...(c_R L\alpha)^{q_R} \\ 
&& = \sum_{q_1,...,q_R ; q_r \geq 1} (q_1+...+q_R) C^r_{q_1,...,q_R} c_1^{q_1}...c_r^{q_r -1}...c_R^{q_R} (L\alpha)^{q_1+...+q_R-1} \\
&& = \frac{1}{c_r} \sum_{q_1,...,q_R ; q_r \geq 1} (q_1+...+q_R) C^r_{q_1,...,q_R} c_1^{q_1}...c_r^{q_r}...c_R^{q_R} (L\alpha)^{q_1+...+q_R-1} \\
&& = \frac{1}{c_r} \frac{ \partial}{\partial L\alpha}  G^r_{LL}(L,\alpha).
\end{eqnarray}
Hence, we have to solve the following system of differential equations:
\begin{equation}
\frac{ \partial}{\partial L\alpha}  G^r_{LL}(L,\alpha) = c_r G_{LL}^r \prod_{r'=1}^R (G_{LL}^{r'})^{\eta_{r'}}.
\end{equation}
This is quite trivial and is achieved defining $x=\alpha L$ and $y_r(x) = \ln G_{LL}^r(\alpha,L)$. The boundary conditions are $y^r(0) = 0$.  Then we can write
\begin{eqnarray}
y'_r &=& c_r \exp\left( \sum_{r'=1}^R \eta_{r'}y_{r'} \right)\\
y'_{r_1} &=& \frac{c_{r_1}}{c_{r_2}} y'_{r_2} \\
y_{r_1} &=& \frac{c_{r_1}}{c_{r_2}} y_{r_2} \\
y'_r &=& c_r \exp\left( \sum_{r'=1}^R \eta_{r'} \frac{c_{r'}}{c_r}y_{r} \right) \\
A &=& \sum_{r'=1}^R \eta_{r'} c_{r'} \\
y'_r &=& c_r \exp\left(  \frac{A}{c_r}y_{r} \right)
\end{eqnarray}
We solve the last differential equation for y, and finally take the exponential of it to find
\begin{eqnarray}
\nonumber G_{LL}^r = (1 + A \alpha L)^{-\frac{c_r}{A}} \text{  if } \eta_r < 0 \\
\nonumber G_{LL}^r = (1 - A \alpha L)^{\frac{c_r}{A}} \text{  if } \eta_r > 0
\end{eqnarray}

\section{Outlook}

We gave here an example on how one can treat a system of combinatorial Dyson Schwinger equations, using techniques of [7]. This allows to compute the leading log expansion of such a system with a nice system of order 1 differential equations. Further work about it will to consider cases with two scales Green's functions (case of general three points Green's functions) or combining several interactions together, which will also include to consider several coupling constants. The equations one will obtain (if one will !) will be partial differential equations, much more complicated to solves, but this could teach more on the combinatoric structures and approximations of Green's functions.

\section*{Appendix}

Here, we compute a general system of combinatorial Dyson Schwinger equation of the form
\begin{equation}
X^r = \mathds{1}_H + \alpha B_+^{\gamma^r}[f^r(X^1,...,X^R)],
\end{equation}
with $f^r$ a function which can be expanded in a formal series in each of its variable. We proves with a (quite tedious) direct derivation that the leading log expansion of the solution of this system is a solution of the following system:
\begin{eqnarray}
\nonumber \sum_{q_1,...,q_R ; q_r \geq 1} (q_1+...+q_R) C^r_{q_1,...,q_R} \phi_R(a_1)^{q_1}...\phi_R(a_r)^{q_r -1}...\phi_R(a_R)^{q_R} \alpha^{q_1+...+q_R} \\ 
 = f^r(G^1_{LL},...,G^R_{LL}).
\end{eqnarray}
About the notation, we write $\delta (m - n)$ for the Kronecker delta $\delta_{m,n}$. We also asserts that Proposition 10, Corollary 3 and Proposition 11 are still truth in this context; in fact, the preceding proofs can be easily extended to this case.

\begin{eqnarray}
\nonumber && X^r = \mathds{1} + \alpha B_+^{a_r}[f^r(X^1,...,X^R)]. \\
\nonumber && \sum_m w_m^r \alpha ^m = \mathds{1} + \alpha B_+^{a_r}\left[\sum_{n_1,...,n_R} f_{n_1,...,n_R}^r {X^1}^{\shuffle_{\Theta}n_1} \shuffle_{\Theta}... \shuffle_{\Theta} {X^R}^{\shuffle_{\Theta} n_R}\right] \\ 
\nonumber && \sum_m w_m^r \alpha ^m = \mathds{1} + \alpha B_+^{a_r} [ \sum_{n_1,...,n_R} f_{n_1,...,n_R}^r \left(\sum_{m_1} w_{m_1}^{r_1}\alpha ^{m_1} \right) ^{\shuffle_{\Theta}n_1} \shuffle_{\Theta}...\\
\nonumber && \shuffle_{\Theta} \left(\sum_{m_R} w_{m_R}^{r_R}\alpha ^{m_R} \right) ^{\shuffle_{\Theta}n_R} ]. \\
\nonumber && \sum_m w_m^r \alpha ^m = \mathds{1} + \alpha B_+^{a_r} [ \sum_{n_1,...,n_R} f_{n_1,...,n_R}^r \left(\sum_{m_{1,1}} w_{m_{1,1}}^{r_1}\alpha ^{m_{1,1}} \right)  \shuffle_{\Theta}... \\
\nonumber && \shuffle_{\Theta} \left(\sum_{m_{1,n_1}} w_{m_{1,n_1}}^{r_1}\alpha ^{m_{1,n_1}} \right) \shuffle_{\Theta}... \shuffle_{\Theta} \left(\sum_{m_{R,1}} w_{m_{R,1}}^{r_R}\alpha ^{m_{R,1}} \right) \shuffle_{\Theta}...\\
\nonumber &&  \shuffle_{\Theta} \left(\sum_{m_{R,n_R}} w_{m_{R,n_R}}^{r_R}\alpha ^{m_{R,n_R}} \right) ].
\end{eqnarray}
Now we can use the linearity of the $B_+$ operator to write:
\begin{eqnarray}
\nonumber && \sum_m w_m^r \alpha ^m = \mathds{1} + \sum_{n_1,...,n_R} f_{n_1,...,n_R}^r \sum_{m_{1,1}}... \\
\nonumber &&\sum_{m_{R,n_R}} \alpha^{1+m_{1,1}+...+m_{R,n_R}} B_+^{a_r} \left[  w_{m_{1,n_1}}^{r_1}  \shuffle_{\Theta}... \shuffle_{\Theta} w_{m_{R,n_R}}^{r_R} \right].
\end{eqnarray}
This allows us to take the equality order by order
\begin{eqnarray}
\nonumber && w^r_0 = \mathds{1}, \\
\nonumber && w^r_m = \sum_{n_1,...,n_R} f_{n_1,...,n_R} \sum_{m_{1,1}}...\sum_{m_{R,n_R}} \delta (1+m_{1,1}+...+m_{R,n_R}-m) \\
\nonumber && B_+^{a_r} \left[  w_{m_{1,1}}^{r_1}  \shuffle_{\Theta}... \shuffle_{\Theta} w_{m_{R,n_R}}^{r_R} \right].
\end{eqnarray}
From Theorem 1, one knows that one can consider only words written as $a_1^{\shuffle_{\Theta} q_1}\shuffle_{\Theta}...\shuffle_{\Theta}a_R^{\shuffle_{\Theta} q_R}$ in $w_m^r$: 
\begin{eqnarray}
\nonumber w_m^r = \sum_{q_1,...,q_R} C^r_{q_1,...,q_R} a_1^{\shuffle_{\Theta} q_1}\shuffle_{\Theta}...\shuffle_{\Theta}a_R^{\shuffle_{\Theta} q_R} \delta(q_1 + ... + q_R - m) + O(m-1),
\end{eqnarray} 
where $O(m-1)$ stands for words of kinematical degree at most $m-1$. One can then rewrite, keeping only the words of kinematical degree m, 
\begin{eqnarray}
\nonumber && \sum_{q_1,...,q_R} C^r_{q_1,...,q_R} a_1^{\shuffle_{\Theta} q_1}\shuffle_{\Theta}...\shuffle_{\Theta}a_R^{\shuffle_{\Theta} q_R} \delta(q_1 + ... + q_R - m) = \\
\nonumber && \sum_{n_1,...,n_R} f_{n_1,...,n_R}^r \sum_{m_{1,1}}...\sum_{m_{R,n_R}} \delta (1+m_{1,1}+...+m_{R,n_R}-m)\\ 
\nonumber && B_+^{a_r} [  \sum_{q_{1,m_{1,1}},...,q_{R,m_{1,1}}} C^{r_1}_{q_{1,m_{1,1}},...,q_{R,m_{1,1}}} a_1^{\shuffle_{\Theta} q_{1,m_{1,1}}} \delta(q_{1,m_{1,1}}+...+q_{R,m_{1,1}} - m_{1,1}) \shuffle_{\Theta}... \\
\nonumber && \shuffle_{\Theta}a_R^{\shuffle_{\Theta} q_{R,m_{1,1}}}  \shuffle_{\Theta}... \shuffle_{\Theta} \sum_{q_{1,m_{R,n_R}},...,q_{R,m_{R,n_R}}} C^{r_R}_{q_{1,m_{R,n_R}},...,q_{R,m_{R,n_R}}}\\
\nonumber &&  a_1^{\shuffle_{\Theta} q_{1,m_{R,n_R}}}\shuffle_{\Theta}...\shuffle_{\Theta}a_R^{\shuffle_{\Theta} q_{R,m_{R,n_R}}} \delta(q_{1,m_{R,n_R}}+...+q_{R,m_{R,n_R}} - m_{R,n_R})].
\end{eqnarray}
One more time one uses the linearity of the $B_+$ operator to write
\begin{eqnarray}
\nonumber && \sum_{q_1,...,q_R} C^r_{q_1,...,q_R} a_1^{\shuffle_{\Theta} q_1}\shuffle_{\Theta}...\shuffle_{\Theta}a_R^{\shuffle_{\Theta} q_R} \delta(q_1 + ... + q_R - m) = \\
\nonumber && \sum_{n_1,...,n_R} f_{n_1,...,n_R}^r \sum_{m_{1,1}}...\sum_{m_{R,n_R}} \delta (1+m_{1,1}+...+m_{R,n_R}-m)\\ 
\nonumber && \sum_{q_{1,m_{1,1}},...,q_{R,m_{1,1}}} C^{r_1}_{q_{1,m_{1,1}},...,q_{R,m_{1,1}}} \delta(q_{1,m_{1,1}}+...+q_{R,m_{1,1}} - m_{1,1})...  \\
\nonumber &&  \sum_{q_{1,m_{R,n_R}},...,q_{R,m_{R,n_R}}} C^{r_R}_{q_{1,m_{R,n_R}},...,q_{R,m_{R,n_R}}} \delta(q_{1,m_{R,n_R}}+...+q_{R,m_{R,n_R}} - m_{R,n_R})\\
\nonumber &&  B_+^{a_r} [ a_1^{\shuffle_{\Theta} q_{1,m_{1,1}}+...+ q_{1,m_{R,n_R}}}\shuffle_{\Theta}...\shuffle_{\Theta}a_R^{\shuffle_{\Theta} q_{R,m_{1,1}}+...+ q_{R,m_{R,n_R}}}  ].
\end{eqnarray}
Use Proposition 11 to rewrite:
\begin{eqnarray}
\nonumber && (q_1+...+q_R) C^r_{q_1,...,q_R} a_1^{\shuffle_{\Theta} q_1}\shuffle_{\Theta}...\shuffle_{\Theta}a_R^{\shuffle_{\Theta} q_R} \delta(q_1 + ... + q_R - m) = \\
\nonumber &&\sum_{n_1,...,n_R} f_{n_1,...,n_R}^r \sum_{m_{1,1}}...\sum_{m_{R,n_R}} \delta (1+m_{1,1}+...+m_{R,n_R}-m)\\ 
\nonumber && \sum_{q_{1,m_{1,1}},...,q_{R,m_{1,1}}} C^{r_1}_{q_{1,m_{1,1}},...,q_{R,m_{1,1}}} \delta(q_{1,m_{1,1}}+...+q_{R,m_{1,1}} - m_{1,1})...  \\
\nonumber && \sum_{q_{1,m_{R,n_R}},...,q_{R,m_{R,n_R}}} C^{r_R}_{q_{1,m_{R,n_R}},...,q_{R,m_{R,n_R}}} \delta(q_{1,m_{R,n_R}}+...+q_{R,m_{R,n_R}} - m_{R,n_R}) \\
\nonumber &&a_1^{\shuffle_{\Theta} q_{1,m_{1,1}}+...+ q_{1,m_{R,n_R}}}\shuffle_{\Theta}... \shuffle_{\Theta} a_r^{\shuffle_{\Theta} q_{r,m_{1,1}}+...+ q_{r,m_{R,n_R}} +1 } \shuffle_{\Theta}...\\
\nonumber && \shuffle_{\Theta}a_R^{\shuffle_{\Theta} q_{R,m_{1,1}}+...+ q_{R,m_{R,n_R}}} \delta(q_{1,m_{1,1}}+...+ q_{1,m_{R,n_R}} -q_1)... \\
\nonumber && \delta(q_{r,m_{1,1}}+...+ q_{r,m_{R,n_R}} +1 -q_r)... \delta(q_{R,m_{1,1}}+...+ q_{R,m_{R,n_R}}-q_R).
\end{eqnarray}
The next step is to apply the Feynman rules. Moreover, we reorganize the sum in order to put each coefficient $C$ in front of its corresponding term in the sum.
\begin{eqnarray}
\nonumber && (q_1+...+q_R) C^r_{q_1,...,q_R} \phi_R(a_1)^{q_1}...\phi_R(a_R)^{q_R} \delta(q_1 + ... + q_R - m) = \\
\nonumber && \sum_{n_1,...,n_R} f_{n_1,...,n_R}^r \sum_{m_{1,1}}...\sum_{m_{R,n_R}} \delta (1+m_{1,1}+...+m_{R,n_R}-m)\\ 
\nonumber && \sum_{q_{1,m_{1,1}},...,q_{R,m_{1,1}}} C^{r_1}_{q_{1,m_{1,1}},...,q_{R,m_{1,1}}} \phi_R(a_1)^{q_{1,m_{1,1}}}...\phi_R(a_R)^{q_{R,m_{1,1}}}\\
\nonumber && \delta(q_{1,m_{1,1}}+...+q_{R,m_{1,1}} - m_{1,1})...  \\
\nonumber && \sum_{q_{1,m_{R,n_R}},...,q_{R,m_{R,n_R}}} C^{r_R}_{q_{1,m_{R,n_R}},...,q_{R,m_{R,n_R}}} \phi_R(a_1)^{q_{1,m_{R,n_R}}}...\phi_R(a_R)^{q_{R,m_{R,n_R}}}\\
\nonumber && \delta(q_{1,m_{R,n_R}}+...+q_{R,m_{R,n_R}} - m_{R,n_R}) \phi_R(a_r)  \delta(q_{1,m_{1,1}}+...+ q_{1,m_{R,n_R}} -q_1)...\\
\nonumber &&  \delta(q_{r,m_{1,1}}+...+ q_{r,m_{R,n_R}} +1 -q_r)... \delta(q_{R,m_{1,1}}+...+ q_{R,m_{R,n_R}}-q_R).
\end{eqnarray}
Then one divide by $\phi_R(a_r)$ both sides, which is allowed for $n\geq 1$ (in the trivial case $n=0$, $\phi_R(w_0^r) = 1$):
\begin{eqnarray}
\nonumber && (q_1+...+q_R) C^r_{q_1,...,q_R} \phi_R(a_1)^{q_1}...\phi_R(a_r)^{q_r -1}...\phi_R(a_R)^{q_R} \delta(q_1 + ... + q_R - m) =\\
\nonumber && \sum_{n_1,...,n_R} f_{n_1,...,n_R}^r \sum_{m_{1,1}}...\sum_{m_{R,n_R}} \delta (1+m_{1,1}+...+m_{R,n_R}-m) \\ 
\nonumber && \sum_{q_{1,m_{1,1}},...,q_{R,m_{1,1}}} C^{r_1}_{q_{1,m_{1,1}},...,q_{R,m_{1,1}}} \phi_R(a_1)^{q_{1,m_{1,1}}}...\phi_R(a_R)^{q_{R,m_{1,1}}} \\
\nonumber && \delta(q_{1,m_{1,1}}+...+q_{R,m_{1,1}} - m_{1,1})...\\
\nonumber && \sum_{q_{1,m_{R,n_R}},...,q_{R,m_{R,n_R}}} C^{r_R}_{q_{1,m_{R,n_R}},...,q_{R,m_{R,n_R}}} \phi_R(a_1)^{q_{1,m_{R,n_R}}}...\phi_R(a_R)^{q_{R,m_{R,n_R}}}\\\nonumber &&\delta(q_{1,m_{R,n_R}}+...+q_{R,m_{R,n_R}} - m_{R,n_R})  \\
\nonumber &&  \delta(q_{1,m_{1,1}}+...+ q_{1,m_{R,n_R}} -q_1)... \delta(q_{r,m_{1,1}}+...+ q_{r,m_{R,n_R}} +1 -q_r)... \\
\nonumber && \delta(q_{R,m_{1,1}}+...+ q_{R,m_{R,n_R}}-q_R).
\end{eqnarray}
And we do the whole resummation over $q_1,...q_R$, re-injecting also the $\alpha$ ; the problematic case $q_r=0$ being cancelled by $\delta(q_{r,m_{1,1}}+...+ q_{r,m_{R,n_R}} +1 -q_r)$. Hence,
\begin{eqnarray}
\nonumber \sum_{q_1,...,q_R ; q_r \geq 1} (q_1+...+q_R) C^r_{q_1,...,q_R} \phi_R(a_1)^{q_1}...\phi_R(a_r)^{q_r -1}...\phi_R(a_R)^{q_R} \alpha^{q_1+...+q_R} \\ 
\nonumber = f^r(G^1_{LL},...,G^R_{LL}).
\end{eqnarray}

\section*{Bibliography}

\vspace{5pt}

[1] Christophe Reutenauer, \textit{Free Lie algebras}, Oxford University Press (1993)

\vspace{5pt}

[2] Michael E. Hoffman, \textit{Quasi shuffle product}, Journal of Algebraic Combinatorics \textbf{11} (2000), 49-68

\vspace{5pt}

[3] Dirk Kreimer, \textit{Dyson Schwinger equations : From Hopf algebra to number theory}, Universality and renormalization, vol \textbf{50} of Fields Institute Communications, pp. 225-248,  Am. Math. Soc., Providence (2007)

\vspace{5pt}

[4] Dirk Kreimer, Lutz Klaczynski \textit{Dyson Schwinger Equations, Dirk Kreimer's Lecture, notes by Lutz Klaczynski}, www2.mathematik.hu-berlin.de/~kreimer/wp-content/uploads/SkriptRGE.pdf

\vspace{5pt}

[5] Joachim Kock \textit{Polynomial functors and Combinatorial Dyson Schwinger equations}, Preprint, arXiv:1512.03027

\vspace{5pt}

[6] Loic Foissy \textit{General Dyson Schwinger Equations and system},  Comm. Math. Phys. \textbf{327} (2014), 151-179

\vspace{5pt}

[7] Olaf Krueger, Dirk Kreimer \textit{Filtrations in Dyson Schwinger equations : next-to$^{\lbrace j \rbrace}$-leading log expansions systematically}, Annals of Physics (2015)

\vspace{5pt}

[8] Dominique Manchon \textit{Hopf algebras, from basics to applications to renormalization}, Comptes-rendus des Rencontres mathématiques de Glanon de 2001 (2006)

\vspace{5pt}

[9] Alain Connes, Dirk Kreimer \textit{Hopf Algebras, Renormalization and Noncommutative Geometry}, Comm. Math. Phys. \textbf{199} (1998), 203-242

\end{document}